\documentclass[12pt,oneside]{article}
\usepackage{amsmath,amssymb,amsfonts,amsthm}
\def\goth{\mathfrak}

\newcommand{\prof}{\noindent \textit{\textbf{Proof.\:\:\,}}}

\def\G{$\Gamma$\,\,}

\def\x{\in{\goth X}(M)}

\def\Section#1{\vspace{30truept}\addtocounter{section}{1}\setcounter{thm}{0}
\setcounter{equation}{0}{\noindent\Large\bf
    \arabic{section}.~~#1}\par \vspace{12pt}}

\newtheorem{thm}{Theorem}[section]
\newtheorem{cor}[thm]{Corollary}
\newtheorem{lem}[thm]{Lemma}
\newtheorem{prop}[thm]{Proposition}
\newtheorem{defn}[thm]{Definition}

\newtheorem{rem}[thm]{Remark}

\numberwithin{equation}{section}

\title{L-REGULAR LINEAR CONNECTIONS}
\author{Nabil L. Youssef  \and Aly A. Tamim}
\date{}

\begin{document}               
\bibliographystyle{plain}
\maketitle                     
\vspace{-1.15cm}
\begin{center}
{Department of Mathematics, Faculty of Science,\\ Cairo University,
Giza, Egypt.}
\end{center}
\vspace{-0.8cm}
\begin{center}
nyoussef@frcu.eun.eg
\end{center}

\vspace{1cm} \maketitle
\smallskip


\vspace{30truept}\centerline{\Large\bf{Introduction}}\vspace{12pt}
\par
An adequate and interesting approach to the theory of nonlinear
connections has been accomplished by Grifone~\cite{3}. His
definition of a nonlinear connection is based on the geometry of the
tangent bundle $T(M)$ of a differentiable manifold $M$. In his
theory, the natural almost-tangent structure $J$ on $T(M)$ (\cite{5}
and~\cite{8}) plays an extremely important role.
\par Anona~\cite{1} generalized the notion of the natural almost-tangent
structure by considering a vector $1$-form $L$ on the manifold
$M$---not on $T(M)$---satisfying certain conditions. As a by-product
of his work, a generalization of some of Grifone's results was
obtained.
\par The first author of the prsent paper, adopting the point of
view of Anona, generalized Grifone's theory of nonlinear
connections~\cite{10}. Grifone's theory can be retrieved
from~\cite{10} by letting $M$ be the tangent bundle of a
differentiable manifold and $L$ the natural almost-tangent structure
$J$ on $M$.
\par In this paper, we still adopt the point of view of Anona and
continue developing the approach established in~\cite{10}. After the
notations and preliminaries (\S 1), the first part (\S 2) of the
work is devoted to the problem of associating to each $L$-regular
linear connection on $M$ a nonlinear $L$-connection on $M$. The
route we have followed is significantly different from that of
Grifone~\cite{4}. Following Tamnou~\cite{8}, we introduce an
almost-complex and an almost-product structures on $M$ by means of a
given $L$-regular linear connection on $M$. The product of these two
structures defines a nonlinear $L$-connection on $M$, which
generalizes Grifone's nonlinear connection~\cite{4}.
\par The seconed part (\S 3) is devoted to the converse problem:
associating to each nonlinear $L$-connection \G on $M$ an
$L$-regular linear connection on $M$; called the $L$-lift of \G. The
existence of this lift is established and the fundamental tensors
associated with it are studied.
\par In the third part (\S 4), we investigate the $L$-lift of a
homogeneous $L$-connection \G, called the Berwald $L$-lift of \G.
Then we particularize our study to the $L$-lift of a conservative
$L$-connection. This $L$-lift enjoys some interesting properties. We
finally deduce various identities concerning the curvature tensors
of such a lift. This generalizes similar identities found
in~\cite{9}.

\Section{Notations and Preliminaries} The following notations will
be used throughout the paper:\newline $M$: a differentiable manifold
of class $C^\infty$ and of finite dimension.\newline $T(M)$: the
tangent bundle of $M$.\newline ${\goth X}(M)$: the Lie algebra of
vector fields on $M$.\newline $J$: the natural almost-tangent
structure on $T(M)$ (\cite{8} and~\cite{5}).\newline $i_K$: the
interior product with respect to the vector form $K$.

\par All geometric objects considered in this paper are supposed to
be of class $C^\infty$. The formalism of
Fr\"{o}licher-Nijenhuis~\cite{2} will be our fundamental tool. The
whole work is based on the approach developed in~\cite{10}, which
relies, in turn, on~\cite{1} and~\cite{3}. We give here a brief
account of such approach.

\par Let $M$ be a $C^\infty$ manifold of dimension $2n$. Let $L$
be a vector $1$-form on $M$ of constant rank $n$ and such that
$[L,L]=0$ and that $Im(L_z)=Ker(L_z)$ for all $z\in M$. It follows
that $L^2=0$ and $[C,L]=-L$, where $C$ is the canonical vector field
on $M$ ~\cite{10}. We call the linear space $Im(L_z)=Ker(L_z)$ the
{\it vertical space} of $M$ at $z$ and denote it by $V_z(M)$; and as
a vector bundle, we write $V(M)$.

\par A vector form $K$ on $M$ is said to be {\it homogeneous} of
degree $r$ if $[C,K]=(r-1)K$. It is called {\it $L$-semibasic} if
$LK=0$ and $i_XL=0$ for all $X\in V(M)$.
A vector field $S\x$ is said to be an {\it $L$-semispray} on $M$ if
$LS=C$. An $L$-semispray is an {\it $L$-spray } if it is homogeneous
of degree $2$. The {\it potential} of an $L$-semibasic vector
$k$-form $K$ on $M$ is the $L$-semibasic vector $(k-1)$-form defined
by $K^\circ=i_SK$, where $S$ is an arbitrary $L$-semispray.

\par A vector $1$-form \G on $M$ is called a nonlinear
$L$-connection, or simply an \linebreak{\it $L$-connection}, on $M$
if $L\Gamma=L$ and $\Gamma L=-L$. An $L$-connection \G \,on $M$ is
said to be {\it homogeneous} if it is homogeneous of degree $1$ as a
vector $1$-form. A homogeneous $L$-connection \G\,on $M$ is said to
be {\it conservative} if there exists an $L$-spray $S$ on $M$ such
that $\Gamma=[L,S]$. An $L$-connection \G\,on $M$ defines an
almost-product structure on $M$ such that for all $z\in M$, the
eigenspace of $\Gamma_z$ corresponding to the eigenvalue $(-1)$
coincides with the vertical space $V_z(M)$. The {\it vertical and
horizontal projectors} of \G\,are defined respectively by
$v=\frac12(I-\Gamma)$ and $h=\frac12(I+\Gamma)$ and we thus have the
decomposition $T_z(M)=V_z(M)\oplus H_z(M)$ for all $z\in M$, where
$H_z(M)=Im(h_z)$: the {\it horizontal space} at $z$.

\par Let \G\,be an $L$-connection on $M$. The {\it torsion} of \G
is the $L$-semibasic vector $2$-form $T=\frac12[L,\Gamma]$. The {\it
strong torsion} of \G\,is the $L$-semibasic vector $1$-form
$t=T^\circ+[C,v]$. The strong torsion of \G\,vanishes if, and only
if, \G\,is homogeneous with no torsion. The {\it curvature} of
\G\,is the $L$-semibasic vector $2$-form $\Omega=-\frac12[h,h]$. An
$L$-connection \G on $M$ is {\it strongly flat} if both its
curvature and strong torsion vanish. The vector $1$-form $F$ on $M$
defined by $FL=h$ and $Fh=-L$ defines an almost-complex structure on
$M$ such that $LF=v$. $F$ is called the {\it almost-complex
structure associated with \G}.

\Section{Induced $L$-Connections} In this section we show that an
$L$-regular linear connection $D$ on $M$ induces an $L$-connection
on $M$ and we study such $L$-connection in relation with $D$.
\begin{defn}\label{def.1}{\em\cite{7}}
Let $D$ be a linear connection on $M$. The map
$$K:{\goth X}(M)\longrightarrow{\goth X}(M):X\longmapsto D_XC$$
is called the connection map associated with $D$. That is, $K=DC$.
\end{defn}
\begin{defn}\label{def.2}
A linear connection $D$ on $M$ is said to be $L$-almost-tangent
if\linebreak $DL=0$; that is if
$$D_XLY=LD_XY\quad \forall X,Y\x.$$
\end{defn}
\par For an $L$-almost-tangent connection, $K(X)$ is vertical for
every $X\x$.
\begin{defn}\label{def.3}
A linear connection $D$ on $M$ is said to be $L$-regular if it
satisfies the conditions{\em:}\newline {\em(a)} $D$ is
$L$-slmost-tangent,\newline {\em(b)} the map $V(M)\longrightarrow
V(M):X\longmapsto K(X)$ is an isomorphism on $V(M)$.\newline The
inverse of this map will be denoted by $\varphi$.
\end{defn}
\par For an $L$-regular linear connection, $\varphi\circ K=K\circ
\varphi=I$ on $V(M)$.

\bigskip
\par Let $D$ be an $L$-regular linear connection on $M$. By
definition, the vertical component $vX$ of $X\x$ is
$$ vX=\varphi(K(X))$$
and the horizontal component $hX$ is
$$hX=X-\varphi(K(X))$$
Hence, any vector field $X\x$ can be written as $X=vX+hX$ and we
have the decomposition of $T(M)$:
$$T(M)=V(M)\oplus H(M),$$
where $H(M)$ is the vector bundle of horizontal vectors.
\par The vertical and horizontal projectors $v$ and $h$ are thus
given by:
\begin{equation}\label{eq.1}
v=\varphi\circ K,\qquad h=I-\varphi\circ K
\end{equation}
\par One can easily show that:
\begin{equation}\label{eq.2}
Lv=0,\quad vL=L,\quad Lh=L,\quad hL=0
\end{equation}
\begin{equation}\label{eq.3}
K(vX)=K(X),\quad K(hX)=0\quad \forall X\x
\end{equation}
\begin{lem}\label{le.1}\quad
If ${\bf T}$ and ${\bf R}$ are the torsion and curvature tensors of
an $L$-almost-tangent connection $D$ on $M$, respectively,
then\newline {\em(a)} ${\bf T}(LX,LY)=L{\bf T}(LX,Y)+L{\bf
T}(X,LY)$\newline {\em(b)} ${\bf R}(X,Y)LZ=L{\bf R}(X,Y)Z$\newline
for every $X,Y,Z\x$.
\end{lem}
\prof (a) follows from the fact that $D$ is $L$-almost-tangent and
that $L^2=0=[L,L]$.\newline (b) is a direct consequence of the
$L$-almost-tangency of $D$.\ \ $\Box$

\bigskip
\par Let $D$ be an $L$-regular linear connection on $M$.
Following Tamnou~\cite{8}, we will define on $M$ an almost-complex
and an almost-product structures using $L$ and the horizontal
projector $h$ associated with $D$.
\par
Define the vector $1$-form $G$ on $M$ by
\begin{equation}\label{eq.4}
G(LX)=-hX,\quad G(hX)=LX\qquad \forall X\x
\end{equation}
Clearly, $G^2=-I$ and so $G$ is an almost-complex structure on $M$.
\par Using equations (\ref{eq.2}) and (\ref{eq.4}) together with the properties of
$L$, $v$ and $h$, one \linebreak can prove
\begin{prop}\label{pp.1}
The almost-complex structure $G$ has the following
properties{\em:}\newline \noindent $
\begin{array}{@{(}c@{)\ }l@{\qquad\qquad}@{(}c@{)\ }l}
a&GL=-h,\ Gh= L,     &b& LG=-v,      \\
c&Gv=hG = G-L,       &d& vG=G-Gv = L,\\
e&GL+LG=-I,          &f& Gh+hG=G.
\end{array}
$
\end{prop}
\par Again, define the vector $1$-form $H$ on $M$ by
\begin{equation}\label{eq.5}
H(LX)=hX,\quad H(hX)=LX\qquad \forall X\x
\end{equation}
Cleary, $H^2=I$ and so $H$ is an almost-product structure on $M$.
\par Using equations (\ref{eq.2}) and (\ref{eq.5}) together with
Proposition \ref{pp.1} and the properties of $L$, $v$, $h$ and $G$,
one can prove the analogue of Proposition \ref{pp.1} for $H$:
\begin{prop}\label{pp.2}
The almost-product structure $H$ has the following
properties{\em:}\newline \noindent $
\begin{array}{@{(}c@{)\ }l@{\qquad\qquad}@{(}c@{)\ }l}
a&HL=h,\ Hh= L,          &b& LH=v,       \\
c&Hv=hH = H-L=-hG,       &d& vH=H-Hv = L,\\
e&HL+LH=I,               &f& Hh+hH=H,    \\
g&GH=-HG,                &h& G+H=2L.
\end{array}
$
\end{prop}
\par The above Properties (c), (g) and (h) above relate the two structures
$G$ and $H$.
\par The concept of almost-quaternionian structure in the next
result is taken in the sense of Libermann~\cite{6}.
\begin{prop}\label{pp.3}
The pair $(G,H)$ defines an almost-quaternionian structure
on\linebreak $(M,L,D)$.
\end{prop}
\par In fact, $G^2=-H^2=-I\,$ and $\, GH+HG=0$.

\bigskip
\par Now, we define another almost-product structure $\Gamma$, of
extreme importance, in terms of the two structures $G$ and $H$.
\begin{prop}\label{pp.4}
The vector $1$-form $\Gamma=HG$ is an almost-product structure on
$M$.
\end{prop}
\prof Using Proposition \ref{pp.2} and the fact that $G^2=-I$,
$H^2=I$, we get\newline
$\Gamma^2=(HG)(HG)=H(GH)G=-H(HG)G=-H^2G^2=I.$\ \ $\Box$
\begin{thm}\label{th.1}
To each $L$-regular linear connection $D$ on $M$ there is associated
a unique $L$-connection \G on $M$ given by $\Gamma=HG$, where $G$
and $H$ are defined respectively by {\em(\ref{eq.4})} and
{\em(\ref{eq.5})}.
\end{thm}
\prof Using Propositions \ref{pp.1} and \ref{pp.2}, we get:\newline
$L\Gamma=L(HG)=(LH)G=vG=L$ and $\Gamma L=(HG)L=H(GL)=-Hh=-L.$ Hence,
\G is an $L$-connection on $M$. Uniqueness is straightforward.\ \
$\Box$
\begin{defn}\label{def.4}
Let $D$ be an $L$-regular linear connection on $M$. The
$L$-connection \G defined in theorem \ref{th.1} is said to be the
$L$-connection on $M$ induced by $D$.
\end{defn}
\par The next result expresses \G in an explicit form in terms
of the connection map $K$ of Definition \ref{def.1}.
\begin{thm}\label{th.2}
Let $D$ be an $L$-regular linear connection on $M$. The
$L$-connection \G induced by $D$ is expressed in the form
\begin{equation}\label{eq.6}
\Gamma=I-2\varphi\circ K,
\end{equation}
where $K$ is the connection map associated with $D$ and $\varphi$ is
the inverse map of the restriction of $K$ on $V(M)$.
\end{thm}
\prof Using Propositions \ref{pp.1} and \ref{pp.2}, we get:
$$H+G=2L\Longrightarrow\Gamma-I=2LG\Longrightarrow\Gamma-v-h=-2v
\Longrightarrow\Gamma=h-v.$$ Now, for every $X\x$, $\Gamma
X=hX-vX=X-2\varphi(K(X))=(I-2\varphi\circ K)X$; by virtue of
(\ref{eq.1}). Hence (\ref{eq.6}) holds.\ \ $\Box$
\begin{cor}\label{co.1}
We have\:\newline {\em(a)} $\Gamma=h-v.$\newline {\em(b)} $\Gamma
h=h\Gamma=h,\quad \Gamma v=v \Gamma=-v.$
\end{cor}
\begin{cor}\label{co.2}
The vertical and horizontal projectors of\, \G coincide with the
vertical and horizontal projectors of $D$, respectively.
\end{cor}
\par In fact, we have ,
$\frac12(I-\Gamma)=\frac12(I-I+2\varphi\circ K)=\varphi\circ K=v$,
by (\ref{eq.1}) and (\ref{eq.6}). Similarly, $\frac12(I+\Gamma)=h$.
\begin{rem}
{\em When $M=T(N)$; $N$ being a differentiable manifold of dimension
$n$, and $L=J$, the induced nonlinear connection on $M$ defined by
Grifone~\cite{4} is retrieved as a special case of Theorem
\ref{th.2}.}
\end{rem}

\medskip
\par Throughout the remaining part of this section, $D$ will denote an
$L$-regular linear connection on $M$, $K$ its connection map and  \G
the $L$-connection on $M$ induced \linebreak by $D$.
\begin{prop}\label{pp.5}
The $L$-connection \G is homogeneous if, and only if, $K$ is
homogeneous of degree one.
\end{prop}
\prof  We have by (\ref{eq.1}),
\begin{equation}\label{eq.7}
[C,v]=[C,\varphi\circ K]=\varphi\circ[C,K]+[C,\varphi]\circ K
\end{equation}
We calculate the last term of (\ref{eq.7}). For every $X\x$,
$$[C,\varphi]K(X)=[C,(\varphi\circ K)X]-\varphi[C,K(X)]$$
Since $K(X)$ and $[C,(\varphi\circ K)X]$ are vertical and since
$\varphi\circ K=K\circ\varphi=I$ on the vertical bundle, then
$$[C,\varphi]K(X)=(\varphi\circ K)[C,(\varphi\circ K)X]-
\varphi[C,(K\circ\varphi)K(X)]=-\varphi([C,K](\varphi\circ K)X).$$
Hence, we obtain
\begin{equation}\label{eq.8}
[C,\varphi]\circ K=-\varphi\circ[C,K]\circ\varphi\circ K
\end{equation}
It follows from (\ref{eq.7}) and (\ref{eq.8}) that
$[C,v]=\varphi\circ[C,K] \circ(I-\varphi\circ K)$. Then, by
(\ref{eq.1}),
$$[C,v]=\varphi\circ[C,K]\circ h,$$
from which $\;[C,\Gamma]=0\Longleftrightarrow[C,K]=0$.\ \ $\Box$
\begin{defn}\label{def.5}
The connection $D$ is said to be reducible if $\, D\Gamma=0$.
\end{defn}
\par Clearly, $D\Gamma=0$ if, and only if, $Dh=Dv$
\begin{lem}\label{le.2}
Let $F$ be the almost-complex structure associated with \G. A
sufficient condition for $D$ to be reducible is that $DF=0$.
\end{lem}
\prof Corollary \ref{co.1} and the definition of $F$ are used in the
proof.\newline For every $X,Y\x$, we have
\begin{eqnarray*}
D_X \Gamma Y &=& D_X h\Gamma Y+D_X v\Gamma Y=D_X hY-D_X vY=D_XFLY-D_X vY\\
             &=& FD_X LY-D_X vY, \\
\Gamma D_X Y &=& \Gamma D_X hY+\Gamma D_X vY=\Gamma D_XFLY+\Gamma D_X vY\\
             &=& \Gamma FD_X LY+\Gamma D_X vY= FD_X LY-D_X vY,
\end{eqnarray*}
since $\,FD_XLY\,$ is horizontal and $\,D_XvY\,$ is vertical. Hence
the result.\ \ $\Box$
\par The condition of Lemma \ref{le.2} will be shown later to be
necessary (Proposition \ref{pp.6} below).

\begin{thm}\label{th.3}
Let $\overline{D}$ be an $L$-regular linear connection on the vector
bundle \linebreak $V(M)\longrightarrow M$. There exists a unique
reducible connection $D$ on $M$ whose restriction to $V(M)$
coincides with $\overline{D}$.
\end{thm}
\prof  It should first be noticed that for every $X\x$ the operator
$\overline{D}_X$ acts on vertical vector fields while the operator
$D_X$ (to be determined) acts on vector fields on $M$.
\par Let $\overline{K}=\overline{D}C$ and $\overline{\varphi}$
the inverse of the isomorphism of $V(M)$ defined by the restriction
of $\overline{K}$ to $V(M)$. The vector $1$-form
$\overline{\Gamma}=I-2\overline{\varphi}\circ\overline{K}$ is
clearly an $L$-connection on $M$. Let $F$ denote the almost-complex
structure associated with $\overline{\Gamma}$. Set
\begin{equation}\label{eq.9}
D_XY=F\overline{D}_XLY+\overline{D}_XLFY.
\end{equation}
$D$ is a linear connection on $M$ with the required properties. The
proof follows the same lines as in \cite{4} with the necessary
modifications.
\par It is a simple matter to show that $DC=\overline{D}C$.
Consequently, the $L$-connection \G induced by $D$ coincides with
$\overline{\Gamma}$. (This justifies the use of the same symbol $F$
for both almost-complex structurs associated with \G and
$\overline{\Gamma})$.\ \ $\Box$
\begin{prop}\label{pp.6}
The following assertions are equivalent\:\newline {\em(a)} $D$ is
reducible.\newline {\em(b)} $DF=0$.\newline {\em(c)} $Dv=Dh=0.$
\end{prop}
\prof \newline (a)$\Longrightarrow$(b): follows from formula
(\ref{eq.9}).\newline (b)$\Longrightarrow$(c): $Dv=D(LF)=LDF=0, \ \
Dh=D(FL)=FDL=0$.\newline (c)$\Longrightarrow$(a):
$D\Gamma=D(h-v)=Dh-Dv=0$.\ \ $\Box$
\begin{rem}
{\em If an $L$-regular linear connection on $M$ is reducible, it is
completely determined by its action on the vertical bundle.
\par In fact, $D_XhY=D_XFLY=FD_XLY$.
}
\end{rem}

\Section{$L$-Lifts and $L$-Connections}
 We have seen that each
$L$-regular linear connection on $M$ induces canonically an
$L$-connection on $M$. We shall investigate here the converse
problem.
\begin{defn}\label{def.6}
A linear connection $D$ on $M$ is said to be $L$-normal if it
satisfies the conditions\:\newline \em{(a)} $D$ is
$L$-almost-tangent,\newline {\em(b)} $D_{LX}C=LX$ for all $X\x$.
\end{defn}
\par An $L$-normal linear connection is clearly $L$-regular. In
fact, the map $LX\longmapsto K(LX)$ in Definition \ref{def.3} is the
identity map, and so $\varphi=I_{V(M)}$.
\begin{lem}\label{le.3}
Let $D$ be an $L$-almost-tangent linear connection on $M$ such that
\linebreak $D_CLX=L[C,X]$ for all $X\x$. The connection $D$ is
$L$-normal if, and only if, ${\bf T}(C,LX)=0$, where {\bf T} is the
torsion of $D$.
\end{lem}
\prof  We have:
\begin{eqnarray*}
{\bf T}(C,LX)&=&D_CLX-D_{LX}C-[C,LX]=L[C,X]-[C,LX]-D_{LX}C\\
             &=&[L,C]X-D_{LX}C=LX-D_{LX}C.
\end{eqnarray*}
Hence, $\;D_{LX}C=LX\Longleftrightarrow {\bf T}(C,LX)=0$.\ \ $\Box$
\begin{defn}\label{def.7}~\par
\noindent --  Let $D$ be a given $L$-normal linear connection on
$M$.
  The $L$-connection \G on $M$  induced by $D$ is called the
$L$-projection of $D$.\\
--  Let \G be a given $L$-connection on $M$. An $L$-normal linear
connection $D$ on $M$ whose $L$-projection coincides with \G is
called an $L$-lift of \G. If $D$ is reducible, it is called a
reducible $L$-lift of \G.

\end{defn}
\par The following result shows (roughly) that
 there is associated a reducible $L$-lift to each
$L$-connection on $M$.

\begin{thm}\label{th.4}
Let \G be an $L$-connection on $M$ and let $B$ be an $L$-semibasic
vector $2$-form on $M$ such that $B^\circ+[C,h]=0$. There exists a
unique reducible $L$-lift $D$ of \G  whose torsion satisfies $\;{\bf
T}(LX,Y)=B(X,Y)\;$ for all $X,Y\x$.
\end{thm}
\prof  Set \vspace{-0.4cm}
\begin{equation}\label{eq.10}
D_XY=h[LY,F]X+L[vY,F]X+FB(X,Y)+B(X,FY) \vspace{-0.4cm}
\end{equation}
where $F$ is the almost-complex structure associated with \G and $v$
and $h$ are respectively the vertical and horizontal projectors of
\G. The connection $D$ defined by (\ref{eq.10}) is the required
$L$-lift of \G. The proof is similar to that of Theorem III,32 of
\cite{4}.\ \ $\Box$
\smallskip
\par As $DF=0$, the connection (\ref{eq.10}) is completely determined by
(cf. Corollary \ref{co.2}):
\begin{equation}\label{eq.11}
\left.
\begin{array}{rcl}
D_{LX}LY &=& L[LX,Y]\\
D_{hX}LY &=& v[hX,LY]+B(X,Y)
\end{array}
\right\}
\end{equation}
or, again, by
\begin{equation}\label{eq.12}
D_XLY=L[vX,Y]+v[hX,LY]+B(X,Y)
\end{equation}
\begin{rem}\label{re.3}
{\em If \G is homogeneous, $[C,h]=0$. Hence, there exists, for every
homogeneous $L$-connection, a canonical reducible $L$-lift
characterized by ${\bf T}(LX,Y)=0$ for all $X,Y\x$.
\par This $L$-lift is called the Berwald $L$-lift of \G.
}
\end{rem}
\begin{rem}\label{re.4}
{\em If $M=T(N)$; $N$ being of dimension $n$, and $L=J$, the
reducible $J$-lift of a $J$-connection \G on $T(N)$ is nothing but
the lift of \G introduced by Grifone \cite{4}. If moreover \G is
homogeneous and we choose $B=0$, the reducible $J$-lift of \G
coincides with the linear extension, in the sense of Theorem
\ref{th.3}, of the usual Berwald connection. This justifies the
adopted terminology. }
\end{rem}

\bigskip
\par In the remaining part of the present section, let \G denote an
$L$-connection on $M$ and $D$ its reducible $L$-lift corresponding
to the L-semibasic vector $2$-form $B$. Also, let $T$, $t$, $\Omega$
and $F$ denote  the torsion, strong torsion, curvature and
associated almost-complx structure of \G, respectively. Let {\bf T}
and {\bf R} be  the torsion and curvature tensors of the linear
connection $D$, respectively.
\begin{prop}\label{pp.7}
The torsion {\bf T} of the $L$-lift $D$ of \G is given, for all
$X,Y\x$, by\:
$${\bf T}(X,Y)=(F\circ T+\Omega)(X,Y)+(i_FB)(X,Y)+2FB(X,Y)$$
\end{prop}
\prof For all $X,Y\x$, we have
\begin{equation}\label{eq.13}
{\bf T}(X,Y)={\bf T}(hX,hY)+{\bf T}(hX,LFY)+{\bf T}(LFX,hY),
\end{equation}
since ${\bf T}(vX,vY)=B(FX,vY)=0$; $B$ being $L$-semibasic.
\par Using (\ref{eq.11}) and the properties of the tensors associated
with \G, we get after some calculations:
\begin{equation}\label{eq.14}
{\bf T}(hX,hY)=h^*[F,F](X,Y)+2FB(X,Y), \vspace{-0.2cm}
\end{equation}
\begin{equation}\label{eq.15}
{\bf T}(hX,LFY)=B(X,FY),
\end{equation}
\begin{equation}\label{eq.16}
{\bf T}(LFX,hY)=B(FX,Y),
\end{equation}
where $h^*[F,F](X,Y)=\frac12[F,F](hX,hY)$.
\par Substituting (\ref{eq.14}), (\ref{eq.15}) and (\ref{eq.16}) into (\ref{eq.13}) and taking the fact
that $h^*[F,F]=F\circ T+\Omega$ [10] into account, the result
follows.\ \ $\Box$
\begin{thm}
A necessary and sufficient condition for the existence of a
symmetric $L$-lift of an $L$-connection \G is that \G be strongly
flat.
\end{thm}
\prof  Suppose that there exists an $L$-lift of \G such that ${\bf
T}=0$. Thus we have $0={\bf T}(LX,Y)=B(X,Y)$. Hence, by Proposition
\ref{pp.7}, $F\circ T+\Omega=0$. But since $F\circ T$ has horizontal
values while $\Omega$ has vertical values, then $T=0$ and
$\Omega=0$. Now, as \G is homogeneous ($[C,\Gamma]=2B^\circ=0$) and
$T=0$, it follows from Corollary 2 of \cite{10} that $t=0$. Hence,
\G is strongly flat.
\par Conversely, if \G is strongly flat, then \G is homogeneous
\cite{10} and the Berwald $L$-lift of \G is evidently symmetric (cf.
Remark \ref{re.3} and (\ref{eq.14})).\ \ $\Box$

\bigskip
\par As the $L$-lift $D$ of an $L$-connection \G is reducible,
the curvature tensor {\bf R} of $D$ is completely determined by the
three semibasic tensors:
\begin{eqnarray*}
R(X,Y)Z&=&{\bf R}(hX,hY)LZ\\
P(X,Y)Z&=&{\bf R}(hX,LY)LZ\\
Q(X,Y)Z&=&{\bf R}(LX,LY)LZ
\end{eqnarray*}
\par Using (\ref{eq.11}) and (\ref{eq.12}), the
properties of the tensors associated with \G and the fact that $B$
is $L$-semibasic, we get after long calculations

\begin{prop}\label{pp.8}
The three curvature tensors $R$, $P$ and $Q$ of $D$ are respectively
given, for all $X,Y,Z\x$, by\:\newline
{\em(a)} $R(X,Y)Z = (D_{LZ}\Omega)(X,Y)+(D_{hY}B)(Z,X)-(D_{hX}B)(Z,Y)\\
\phantom{(a).R(X,Y)Z
=}+B(FB(Z,X),Y)-B(FB(Z,Y),X)+B(FT(X,Y),Z)$.\newline {\em(b)}
$P(X,Y)Z = (D_{LY}B)(Z,X)+v[hX,L[LY,Z]]+v[LZ,[hX,LY]]\newline
\phantom{(a).R(X,Y)Z =}-L[LY,F[hX,LZ]]-L[LZ,F[hX,LY]]$.\newline
{\em(c)} $Q(X,Y)Z =0$.
\end{prop}

\Section{Berwald $L$-Lifts of Homogeneous} \vspace{-0.5cm}
\begin{center}{\Large\bf $L$-Connections}\end{center}
\par In this section, \G will denote a {\bf homogeneous}
$L$-connection on $M$. The reducible $L$-lift of \G characterized by
${\bf T}(LX,Y)=0$, for all $X,Y\x$, is called the Berwald $L$-lift
of \G (cf. Remark \ref{re.3}).
\par By  virtue of (\ref{eq.11}), the Berwald $L$-lift $D$ is completely determined
by:
\begin{equation}\label{eq.17}
\left.
\begin{array}{rcl}
D_{LX}LY &=& L[LX,Y]\\
D_{hX}LY &=& v[hX,LY]
\end{array}
\right\}
\end{equation}
or, again, by
\begin{equation}\label{eq.18}
D_XLY=L[vX,Y]+v[hX,LY]
\end{equation}
\begin{lem}\label{le.4}
The Berwald $L$-lift $D$ of \G is such that\:\newline {\em(a)}
$D_CLX=L[C,X]$.\newline {\em(b)} $[C,DLX]=D[C,LX]$.
\end{lem}
\prof (a) follows from the first formula of (\ref{eq.17}) by letting
$X=S$; an arbitrary $L$-semispray.\newline (b) follows from
(\ref{eq.18}), the properties of $L$ and those of the tensors
associated with \G and from the Jacobi identity.\ \ $\Box$
\begin{rem}\label{re.5}
{\em In view of the above lemma, as the Berwald $L$-lift $D$ of \G
is reducible, $D$ is an "extended connection of directions" in the
sense of Grifone \cite{4} (where $M=T(N)$ and $L=J$). }
\end{rem}
\begin{prop}\label{pp.9}
The torsion tensor of the Berwald $L$-lift of \G is given by\:
$${\bf T}=F\circ T+\Omega$$
\end{prop}
\par This result follows directly from Proposition \ref{pp.7}.
\begin{cor}\label{co.3}
If \G is a conservative $L$-connection on $M$, then
$${\bf T}=\Omega$$
\end{cor}
\par In fact, $T=0$ for conservative $L$-connections \cite{10}.
\begin{prop}\label{pp.10}
The first curvature tensor of the Berwald $L$-lift $D$ of \G
is\linebreak given by\:
\begin{equation}\label{eq.19}
R(X,Y)Z=(D_{LZ}\Omega)(X,Y)
\end{equation}
\end{prop}
\par This result follows immediately from Proposition \ref{pp.8}.
\begin{thm}\label{th.6}
For the Berwald $L$-lift $D$ of $\Gamma$, we have\:
\begin{equation}\label{eq.20}
R(X,Y)S=\Omega(X,Y),
\end{equation}
where $S$ is an arbitrary $L$-semispray on $M$.
\par Consequently, $R=0$ if, and only if, $\Omega=0$.
\end{thm}
\prof  Setting $Z=S$ in (\ref{eq.19}), taking the fact that $\Omega$
is $L$-semibasic into account,\linebreak  we get
$$R(X,Y)S=(D_C\Omega)(X,Y)=D_C\Omega(X,Y)-\Omega(D_ChX,Y)-\Omega(X,D_ChY)$$
Using Lemma \ref{le.4}(a) together with (\ref{eq.18}), we get
\begin{eqnarray*}
R(X,Y)S &=& L[C,F\Omega(X,Y)]-\Omega([C,X],Y)-\Omega(X,[C,Y])\\
        &=& -[C,L]F\Omega(X,Y)+[C,\Omega(X,Y)]-\Omega([C,X],Y)-\Omega(X,[C,Y])\\
        &=& LF\Omega(X,Y)+[C,\Omega](X,Y)\\
        &=& \Omega(X,Y);\mbox{ $\Omega$ being homogeneous of degree $1$
(since \G is).}
\end{eqnarray*}
\par Now, if $\Omega=0$, then $R=0$, by (\ref{eq.19}). Conversely, if
$R=0$, then $\Omega=0$, by (\ref{eq.20}).\newline (Note that we have
shown, in the course of the proof, that $D_C\Omega=\Omega$.)\ \
$\Box$

\bigskip
\par For the rest of the paper, we consider the Berwald $L$-lift
$D$ of a {\bf conservative} $L$-connection \G on $M$.
\par As a conservative $L$-connection \G on $M$ is homogeneous
with no torsion and is of the form $\Gamma=[L,S]$, we may combine
Theorem \ref{th.6} and Theorems 6, 7 and 9 of \cite{10} to obtain
the following result:
\begin{thm}
For the Berwald $L$-lift of a conservative $L$-connection on $M$,
the \linebreak  following assertions are equivalent\:\newline
{\em(a)} $\Omega^\circ=0$.\newline {\em(b)} $\Omega=0$.\newline
{\em(c)} $R=0$.\newline {\em(d)} $[F,F]=0$.\newline {\em(e)} the
horizontal distribution $z\longmapsto H_z(M)$ is completely
integrable.
\end{thm}

\par As for all linear connections, the (classical) Bianchi's
identities for $D$ are \linebreak  given by:
$${\goth S}\; {\bf R}(X,Y)Z={\goth S}\;\{{\bf T}({\bf T}(X,Y),Z)
+(D_X{\bf T})(Y,Z)\},$$
$${\goth S}\;\{ {\bf R}({\bf T}(X,Y),Z)+(D_X{\bf R})(Y,Z)\}=0,$$
where $\goth S$ denotes the cyclic permutation of the vector fields
$X$, $Y$ and $Z$.
\par But since \G is conservative, we have, by Corollary \ref{co.3},
${\bf T}=\Omega$, which is $L$-semibasic. Thus the above identities
reduce to:
\begin{equation}\label{eq.21}
{\goth S}\;{\bf R}(X,Y)Z={\goth S}\;(D_X\Omega)(Y,Z)
\end{equation}
\begin{equation}\label{eq.22}
{\goth S}\;\{ {\bf R}(\Omega(X,Y),Z)+(D_X{\bf R})(Y,Z)\}=0
\end{equation}
\par These two identities give rise to the following useful
identities.
\begin{prop}\label{pp.11}
For the Berwald $L$-lift of a conservative $L$-connection on $M$, we
have for all $X,Y,Z\x${\em:}\newline {\em(a)} ${\goth
S}\;R(X,Y)Z=0$.\newline {\em(b)} ${\goth S}\;(D_{hX}R)(Y,Z)={\goth
S}\;P(X,F\Omega(Y,Z))$.\newline {\em(c)}
$(D_{LZ}R)(X,Y)=(D_{hY}P)(X,Z)-(D_{hX}P)(Y,Z)$.\newline {\em(d)}
$(D_{LZ}P)(X,Y)=(D_{LY}P)(X,Z)$.\newline {\em(e)}
$P(X,Y)Z=P(Y,X)Z=P(Z,X)Y$. ($P$ is symmetric in its three
variables.)
\end{prop}
{\noindent \bf Sketch of the Proof.}\newline (a) Compute
(\ref{eq.21}) for $hX$, $hY$, $hZ$.\newline (b) Compute
(\ref{eq.22}) for $hX$, $hY$, $hZ$.\newline (c) Compute
(\ref{eq.22}) for $hX$, $hY$, $LZ$.\newline (d) Compute
(\ref{eq.22}) for $hX$, $LY$, $LZ$.\newline (e) Compute
(\ref{eq.21}) for $hX$, $hY$, $LZ$.\newline The calculations are too
long but not difficult. So, we omit them.\ \ $\Box$
\begin{cor}\label{co.4}
We have\:\newline {\em(a)} ${\goth
S}\;(D_{hX}\Omega)(Y,Z)=0$.\newline {\em(b)} ${\goth
S}\;(D_{LX}\Omega)(Y,Z)=0$.\newline {\em(c)} ${\goth
S}\;(D_{LX}R)(Y,Z)=0$.
\end{cor}
{\bf Sketch of the Proof.}\newline (a) follows from Proposition
\ref{pp.11}(a).\newline (b) follows from Proposition \ref{pp.11}(a)
and from (\ref{eq.19}).\newline (c) follows from Proposition
\ref{pp.11}(c) and (e).\ \ $\Box$
\begin{rem}
{\em The identities in Proposition \ref{pp.11} and Corollary
\ref{co.4} are similar to those found in \cite{9}. Nevertheless, the
context here is more general and the scope of validity is much
wider. In fact, the above identities are valid for the large class
of $L$-lifts of conservative $L$-connections, while the identities
in \cite{9} are valid only for the Berwald connection as a $J$-lift
of the canonical connection associated with a Finsler space. }
\end{rem}


\begin{thebibliography}{10}
\bibitem{1}
M. Anona: $d_L$-cohomologie et vari\'{e}t\'{e}s feuiet\'{e}es.\\
Th\`{e}se de 3e cycle, Universit\'{e} de Grenoble, 1978.
\bibitem{2}
A. Fr\"{o}licher and A. Nijenhuis: Theory of vector-valued
differential forms, I.\\
Proc. Kon. Ned. Akad., A, {\bf 59}(1956), 338--359.
\bibitem{3}
J. Grifone: Structure presque-tangente et connexions, I.\\
Ann. Inst. Fourier, Grenoble, {\bf 22, 1}(1972), 287--334.
\bibitem{4}
J. Grifone: Structure presque-tangente et connexions, II.\\
Ann. Inst. Fourier, Grenoble, {\bf 22, 3}(1972), 291--338.
\bibitem{5}
J. Klein and A. Voutier: Formes ext\'{e}rieures
g\'{e}n\'{e}ratrices de sprays.\\
Ann. Inst. Fourier, Grenoble, {\bf 18,1}(1968), 241--260.
\bibitem{6}
P. Libermann: Sur le probl\`{e}me d'\'{e}quivalence de certaines
structures infinit\'{e}simales.\\
Annali di Matematica, {\bf 36}(1954).
\bibitem{7}
A. A. Tamim: Finsler and Riemannian structures on a manifold.\\
An. Univ. Timi\c{s}oara, Ser. Mat., {\bf 25,2}(1987),91--98.
\bibitem{8}
T. Tamnou: Geom\'{e}trie diff\'{e}rentielle du fibr\'{e}
tangent--Connexion
de type finsl\'{e}rien.\\
Th\`{e}se de 3e cycle, Universit\'{e} de Grenoble, 1969.
\bibitem{9}
N. L. Youssef: Sur les tenseurs de courbure de la connexion de
Berwald et ses \linebreak distributions de nullit\'{e}.\\
Tensor, N. S., {\bf 36}(1982), 275--280.
\bibitem{10}
N. L. Youssef: $L$-connections and associated tensors.\\
Tensor, N. S., Vol. 60(1998), 229-238. (ArXiv No.: math.DG/0605338).
\end{thebibliography}
\end{document}